\documentclass[12pt]{article}
\usepackage{graphics, latexsym}
\usepackage{graphicx}
\usepackage{setspace}
\usepackage{amssymb}
\usepackage{amsmath}
\usepackage{epsfig}
\usepackage{CJK}
\usepackage{verbatim}
\usepackage{amsmath, amssymb}
\usepackage{makecell}
\usepackage{multirow}
\usepackage{graphicx}
\usepackage{subfigure}
\usepackage{epstopdf}

\usepackage{amssymb}
\usepackage{amssymb}
\usepackage{lineno,hyperref}
\usepackage{url}
\usepackage{bm}
\usepackage{overpic}
\usepackage{graphicx}
\usepackage{mathrsfs}
\usepackage{amsfonts}
\usepackage{amssymb}
\usepackage{amsmath,amsthm}
\usepackage{indentfirst}
\usepackage{footmisc}
\usepackage{multirow}
\usepackage{xcolor}
\usepackage{tikz}
\usepackage{bm}
\usepackage{algorithm}
\usepackage{algpseudocode}
\usepackage{amsmath}
\usepackage{colortbl}
\usepackage{amsmath}
\newtheorem{theorem}{Theorem}
\newtheorem{remark}{Remark}
\usepackage{tabularx}
\usepackage{array}

\usetikzlibrary{arrows,shapes,chains}
\usepackage{pgf}
\usepackage{threeparttable}
\usepackage{diagbox}
\usetikzlibrary{arrows, decorations.pathmorphing, backgrounds, positioning, fit, petri, automata}
\usepackage{tikz}
\usepackage{tabu}
\usepackage{float}
\usetikzlibrary{shadows,arrows,positioning,shapes.geometric} 
\pgfdeclarelayer{background}
\pgfdeclarelayer{foreground}
\pgfsetlayers{background,main,foreground}

\usepackage{booktabs}
\usepackage{array, caption, threeparttable}
\usepackage[font=small,labelfont=bf,labelsep=none]{caption}
\newcolumntype{C}[1]{>{\centering}p{#1}}
\usepackage{makecell}

\begin{document}
\title{ A new method for estimating the tail index using truncated sample sequence $^*$}
\author{Fuquan Tang$\,\,\,\,\,\,$  Dong Han$^\star$ \\
Department of Statistics, School of Mathematical Sciences,\\
Shanghai Jiao Tong University, Shanghai, 200240, China}

\maketitle
\begin{abstract}
This article proposes a new method of truncated estimation to estimate the tail index $\alpha$ of the extremely heavy-tailed distribution with infinite mean or variance. We not only present two truncated estimators $\hat{\alpha}$ and $\hat{\alpha}^{\prime}$ for estimating $\alpha$ ($0<\alpha \leq 1$) and $\alpha$ ($1<\alpha \leq 2$) respectively, but also prove their asymptotic statistical properties. The numerical simulation results comparing the six known estimators in estimating error, the Type I Error  and the power of estimator show that the performance of the two new truncated estimators is quite good on the whole.
\end{abstract}
\renewcommand{\thefootnote}{\fnsymbol{footnote}}
\footnotetext{
\newline$^*$Supported by National Natural Science Foundation of China (11531001)
\newline$^\star$ Corresponding author, E-mail: donghan@sjtu.edu.cn }
\textbf{KEYWORDS}: Heavy-tailed distributions, tail index, truncated sample mean, simulation.

\section*{1. Introduction}
\noindent
Heavy-tailed phenomena are widespread in many aspects of our lives, and exist in a variety of disciplines such as physics, meteorology, computer science, biology, and finance. The probabilistic and statistical methods and theories about the heavy-tailed phenomenon have been used to study the magnitude of earthquakes, the diameter of lunar craters on the surface of the moon, the size of interplanetary fragments, and the frequency of words in human languages, and so on [\ref{bb4}, \ref{bb5}, \ref{bb17}, \ref{bb18}].

Geography and hydrology are important scenarios for the study and application of thick-tailed distribution. In 1998, Anderson [\ref{bb1}] discussed heavy tail time series models and provided a periodic ARMA model for Salt River. In 2022, Merz et al.[\ref{bb22}] provided a detailed and coherent review on understanding heavy tails of flood peak distributions, they proposed nine hypotheses on the mechanisms generating heavy-tailed phenomena in flood system. In financial markets, Mandelbrot[\ref{bb26}] presented seminal research on cotton price using the heavy tails distribution theory. In 2013,  Marat et al.[\ref{bb24}] found the emerging exchange markets would be more pronouncedly heavy-tailed and illustrated that heavy-tailed properties did not change obviously during the financial and economic crisis period.

There is a large literature proposing numerous ideas and methods on the estimation of the tail index $\alpha$ of heavy-tailed distribution. The size of $\alpha$ is mainly used to measure the degree of thinness of the tail. The smaller the $\alpha$, the higher the probability of a heavy-tailed event. Since Hill put forward the famous Hill estimator in 1975 [\ref{bb10}], researchers have provided multiple estimation methods for estimating $\alpha$, such as DPR estimator[\ref{bb15}, \ref{bb16}], QQ estimator[\ref{bb14}], the Moment estimator[\ref{bb6}], $L^p$ quantile estimator[\ref{bb19}], the estimators of extreme value index in a censorship framework[\ref{bb2}, \ref{bb3}, \ref{bb8}, \ref{bb20}], t-Hill estimator[\ref{bb11}, \ref{bb12}], IPO estimator[\ref{bb13}] and so on. There are more than 100 tail index estimators have been reviewed by two papers [\ref{bb7}, \ref{bb9}].

It can be seen that nearly all estimators  based on the order statistics of observation samples. Moreover, the estimators based on the order statistics have three characteristics that are not very satisfactory: (\textbf{1}) The calculation of the estimators is relatively complex since the order statistics are not easy to calculate for large  sample size; (\textbf{2}) The mathematical meaning of the estimators for the tail index is not obvious; (\textbf{3}) There is no explicit expression for the rate of strong consistency convergence of the estimators.

In order to make up for the shortcomings of existing estimation methods, we propose  a new truncated estimation method to estimate the tail index $\alpha$ ($0<\alpha \leq 2$) of heavy-tailed distribution with infinite mean or variance. The proposed two estimators $\hat{\alpha}$ for $0<\alpha \leq  1$ and $\hat{\alpha}^{\prime}$ for $1<\alpha \leq  2$, are based on the truncated sample mean and the truncated sample second order moment, respectively, and they are not only relatively easy to calculate, but also their strong consistency convergence rate and the asymptotic normal property can be obtained.

In Section 2, we will present two truncated estimators $\hat{\alpha}$ and $\hat{\alpha}^{\prime}$, and obtain their asymptotic statistical properties. Section 3 compares the two truncated estimators with the six known estimators in estimating error, the type I error  and the power of estimator by numerical simulations. Section 4 provides concluding remarks. The proofs of the three theorems are given in the Appendix.

\section*{2. Two truncated estimators}
\noindent
Due to a random variable $X$ can be written as the summation of positive and negative parts $X=X^+-X^-$, we consider only the nonnegative random variables in the paper.
Let $X_k, k\geq 1,$ be independent and identical distribution (i.i.d.) with extremely heavy-tailed distribution function $F(x)=1-1/x^{\alpha}$ for $x\geq 1$, where the tail index, $\alpha \in (0, \, 2]$, is unknown. We know that when  $\alpha \in (0, \, 1]$, the mean or variance is infinite, and when $\alpha \in (1, \, 2]$, the mean is finite but the variance is infinite.

In this section, we will present two truncated estimators $\hat{\alpha}$ and $\hat{\alpha}'$ to estimate $\alpha$ ($0<\alpha \leq 1$) and  $\alpha$ ($1<\alpha\leq 2$) respectively and prove their asymptotic statistical properties.

To this end, let  $\{b_n\}$ be a positive truncated sequence satisfying $b_n \nearrow \infty$ as $n\to \infty$. Define the truncated random variable $ X_k(b_n):=X_kI(X_k\leq b_n)$, where $I(\cdot)$ is the indicator function. We can get the truncated  mean $\mu_n$, the truncated sample mean $\hat{\mu}_n$,  truncated second order moment $\nu^2_n$ and the truncated sample second order moment  $\hat{\nu}^2_n$ in the following.
\begin{eqnarray}
\mu_n:=\text{E}(X_k(b_n))= \frac{\alpha}{1-\alpha}\Big(b_n^{1-\alpha}-1\Big),\,\,\,\,\,\,\,\hat{\mu}_n:=n^{-1}\sum_{k=1}^nX_k(b_n)
\label{han1}
\end{eqnarray}
for $0<\alpha <1$ and
\begin{eqnarray}
\nu^2_n:=\text{E}(X^2_k(b_n))=\frac{\alpha}{2-\alpha}\Big(b_n^{2-\alpha}-1\Big),\,\,\,\,\,\,\,\hat{\nu}^2_n:=n^{-1} \sum_{k=1}^nX^2_k(b_n)
\label{han2}
\end{eqnarray}
for $1<\alpha <2$. It follows from equation(\ref{han1}) and (\ref{han2}) that
\begin{eqnarray}
\alpha=1-\frac{\log[\frac{1-\alpha}{\alpha}\mu_n +1]}{\ln b_n}
\label{han3}
\end{eqnarray}
for $0<\alpha <1$ and
\begin{eqnarray}
\alpha=2-\frac{\log[\frac{2-\alpha}{\alpha}\nu^2_n +1]}{\ln b_n}
\label{han4}
\end{eqnarray}
for $1<\alpha <2$.

Hence, we can define two  truncated estimators $\hat{\alpha}$ and $\hat{\alpha}'$ which satisfy the following two equations $x=G_1(x)$ for  $0<x<1$ and $y=G_2(y)$ for  $1<y<2$ respectively,  by replacing $\alpha$, $\mu_n$ and $\nu^2_n$ in equation(\ref{han3}) and (\ref{han4}) with  $\hat{\alpha}$, $\hat{\mu}_n,$  $\hat{\alpha}'$ and $\hat{\nu}^2_n$, respectively, that is,
\begin{eqnarray}
\hat{\alpha}=G_1(\hat{\alpha}):=1-\frac{\ln[\frac{1-\hat{\alpha}}{\hat{\alpha}}\hat{\mu}_n +1]}{\ln b_n}
\label{han5}
\end{eqnarray}
for $0<\hat{\alpha} <1$ and
\begin{eqnarray}
\hat{\alpha}'=G_2(\hat{\alpha}'):=2-\frac{\ln[\frac{2-\hat{\alpha}'}{\hat{\alpha}'}\hat{\nu}^2_n +1]}{\ln b_n}
\label{han6}
\end{eqnarray}
for $1<\hat{\alpha}' <2$.

Take $\hat{\alpha}\nearrow 1$ in equation(\ref{han5}) and $\hat{\alpha}'\nearrow 2$ in equation(\ref{han6}), respectively, it follows that $\hat{\alpha}\sim  \hat{\mu}_n/\ln b_n$ and $\hat{\alpha}' \sim  \hat{\nu}^2_n/\ln b_n$. Hence, we can use the following two estimators $\hat{\alpha}$ and $\hat{\alpha}'$ to estimate $\alpha =1$ and $\alpha =2$, respectively.
\begin{eqnarray}
\hat{\alpha}:=\frac{\hat{\mu}_n}{\ln b_n},
\label{han7}
\end{eqnarray}
\begin{eqnarray}
\hat{\alpha}' := \frac{\hat{\nu}^2_n}{\ln b_n}.
\label{han8}
\end{eqnarray}

Since it is difficult to obtain the analytic solutions (estimators) $\hat{\alpha}$ and $\hat{\alpha}'$ to the two equations  $\hat{\alpha}-G_1(\hat{\alpha})=0$  and $\hat{\alpha}'-G_2(\hat{\alpha}')=0$  respectively,  we present two  recursive estimators for $k\geq 1$ in the following
 \begin{eqnarray}
     \hat{\alpha}_k &=& G_1(\hat{\alpha}_{k-1}),\,\,\,\,\,\,\,\, \text{ if }\,\, 0<\hat{\alpha}_0<1
     \label{han9}
 \end{eqnarray}
 for $0<\alpha\leq 1$ and
 \begin{eqnarray}
     \hat{\alpha}'_k &=& G_2(\hat{\alpha}'_{k-1}),\,\,\,\,\,\,\,\, \text{ if }\,\,  1<\hat{\alpha}'_0<2
     \label{han10}
 \end{eqnarray}
for $1<\alpha \leq 2$, where $\hat{\alpha}_0$ and $\hat{\alpha}'_0$ are two constants.

The following theorem shows that the two estimators $\hat{\alpha}$ and $\hat{\alpha}'$ can be  approximated by the two  sequences  of estimators $\{\hat{\alpha}_k\}$ and $\{\hat{\alpha}'_k\}$, respectively.\\

\begin{theorem}{\label{theorem1}}
{\label{theorem1}} Let $0 \leq \beta<1 / 2$ and $b_{n}$ satisfy $b_{n}^{\alpha} \ln b_{n} \leq \alpha n^{1-2 \beta} / \ln n$ for $0<\alpha<2, \alpha \neq 1$ and large n. Then, both the two equations $x-G_{1}(x)=0,0<x<1$, and $y-G_{2}(y)=0$, $1<y<2$, have unique solutions $\hat{\alpha}$ and $\hat{\alpha}^{\prime}$, respectively. If $0<\hat{\alpha}_{0}<\hat{\alpha}$ and $1<\hat{\alpha}_{0}^{\prime}<\hat{\alpha}^{\prime}<2$ (or $0<\hat{\alpha}<\hat{\alpha}_{0}<1$ and $1<\hat{\alpha}^{\prime}<\hat{\alpha}_{0}^{\prime}<2$ ), then $\hat{\alpha}_{k} \nearrow \hat{\alpha}$ and $\hat{\alpha}_{k}^{\prime} \nearrow \hat{\alpha}^{\prime}$ (or $\hat{\alpha}_{k} \searrow \hat{\alpha}$ and $\hat{\alpha}_{k}^{\prime} \searrow \hat{\alpha}^{\prime}$ ) and
\begin{eqnarray}
\left|\hat{\alpha}-\hat{\alpha}_{k}\right| \leq\left(\frac{1}{\hat{\alpha}^{\star}(1-\hat{\alpha}) \ln b_{n}}\right)^{k}\left|\hat{\alpha}-\hat{\alpha}_{0}\right|
\label{han11}
\end{eqnarray}
\begin{eqnarray}
\left|\hat{\alpha}^{\prime}-\hat{\alpha}_{k}^{\prime}\right| \leq\left(\frac{2}{\hat{\alpha}^{\star \prime}\left(2-\hat{\alpha}^{\prime}\right) \ln b_{n}}\right)^{k}\left|\hat{\alpha}^{\prime}-\hat{\alpha}_{0}^{\prime}\right|
\label{han12}
\end{eqnarray}
where $\hat{\alpha}^{\star}=\min \left\{\hat{\alpha}, \hat{\alpha}_{0}\right\}$ and $\hat{\alpha}^{\star \prime}=\min \left\{\hat{\alpha}^{\prime}, \hat{\alpha}_{0}^{\prime}\right\}$.
\end{theorem}

\begin{remark}{\label{remark1}}
Take large $n$ such that $A:=\hat{\alpha}^{\star}(1-\hat{\alpha}) \ln b_{n}>1$ and $B:=\hat{\alpha}^{\star \prime}\left(2-\hat{\alpha}^{\prime}\right) \ln b_{n} / 2>1$. Note that $\left|\hat{\alpha}-\hat{\alpha}_{0}\right| \leq 1$ and $\left|\hat{\alpha}^{\prime}-\hat{\alpha}_{0}^{\prime}\right| \leq 1$. It follows from the two inequalities  (\ref{han11}) and (\ref{han12}) that
$$
\left|\hat{\alpha}-\hat{\alpha}_{k}\right| \leq e^{-k \ln A}, \quad\left|\hat{\alpha}^{\prime}-\hat{\alpha}_{k}^{\prime}\right| \leq e^{-k \ln B} .
$$
\end{remark}
The two inequalities above implies that $\left\{\hat{\alpha}_{k}\right\}$ and $\left\{\hat{\alpha}_{k}^{\prime}\right\}$ can converge (almost everywhere) at least exponentially to $\hat{\alpha}$ and $\hat{\alpha}^{\prime}$, respectively.

\begin{remark}{\label{remark2}}
If we don't know whether $\alpha$ is included in interval $(0, 1]$ or interval $(1, 2]$, we may take the initial values $\hat{\alpha}_{0}$ and $\hat{\alpha}'_{0}$ in the following: Take  $n_0$ samples (for example, $n_0=50$) such that $\hat{\alpha}_{0}=\frac{\hat{\mu}_{n_0}}{\ln b_{n_0}}$ for $\frac{\hat{\mu}_{n_0}}{\ln b_{n_0}}\leq 1$ and $\hat{\alpha}'_{0}=1.5$ for $\frac{\hat{\mu}_{n_0}}{\ln b_{n_0}}>1$.
\end{remark}

In order to get the asymptotic statistical properties $\hat{\alpha}$ and $\hat{\alpha}^{\prime}$, we  give a theorem in the following, which describes the asymptotic statistical properties of the truncated sample mean $\hat{\mu}_{n}$ and the truncated second order moment $\hat{\nu}_{n}^{2}$.

\begin{theorem}{\label{theorem2}}
\textit{  Assume that the conditions of Theorem 1 hold.    Then
\begin{eqnarray}
   \textbf{P}\Big( \frac{|\hat{\mu}_n-\mu_n|}{\mu_n} \geq \frac{2}{n^{\beta}\sqrt{\ln b_n}}\Big)\leq \frac{2}{n^2},\,\,\,\,\,\,\,\,\,\, \frac{\sqrt{n}(\hat{\mu}_n-\mu_n)}{b_n^{1-\alpha/2}} \Rightarrow N\Big(0,\,\, \frac{\alpha}{2-\alpha}\Big)
\label{han13}
 \end{eqnarray}
for $0<\alpha\leq 1$ and
\begin{eqnarray}
   \textbf{P}\Big( \frac{|\hat{\nu}^2_n-\nu^2_n|}{\nu^2_n} \geq \frac{2}{n^{\beta}\sqrt{\ln b_n}}\Big)\leq \frac{2}{n^2},\,\,\,\,\,\,\,\,\,\, \frac{\sqrt{n}(\hat{\nu}^2_n-\nu^2_n)}{b_n^{2-\alpha/2}} \Rightarrow N\Big(0,\,\, \frac{\alpha}{4-\alpha}\Big)
\label{han14}
 \end{eqnarray}
for $1<\alpha\leq 2$, where "$\Rightarrow$"  denotes the convergence in distribution and $N(\mu, \sigma^2)$ is the normal distribution.}
\end{theorem}
The following theorem gives  the asymptotic statistical properties of the two truncated estimators $\hat{\alpha}$  and $\hat{\alpha}'$.\\

\begin{theorem}{\label{theorem3}}
\textit{ Assume that the conditions of Theorem 1 hold.  Then
\begin{eqnarray}
   \textbf{P}\Big( |\hat{\alpha}-\alpha| \geq \frac{2}{n^{\beta}\ln b_n\sqrt{\ln b_n}}\Big)\leq \frac{2}{n^2},\,\,\,\,\,\,\,\,\,\, \frac{\sqrt{n}(\hat{\alpha}-\alpha)\ln b_n}{b_n^{\alpha/2}} \Rightarrow N\Big(0,\,\, \frac{(1-\alpha)^2}{\alpha(2-\alpha)}\Big)
\label{han15}
 \end{eqnarray}
for $0<\alpha<1$ and
\begin{eqnarray}
   \textbf{P}\Big( |\hat{\alpha}'-\alpha| \geq \frac{2}{n^{\beta}\ln b_n\sqrt{\ln b_n}}\Big)\leq \frac{2}{n^2},\,\,\,\,\,\,\,\,\,\, \frac{\sqrt{n}(\hat{\alpha}'-\alpha)\ln b_n}{b_n^{\alpha/2}} \Rightarrow N\Big(0,\,\, \frac{(2-\alpha)^2}{\alpha (4-\alpha)}\Big)
\label{han16}
 \end{eqnarray}
for $1<\alpha < 2$. Moreover,
\begin{eqnarray}
   \textbf{P}\Big( |\hat{\alpha}-1| \geq \frac{2\sqrt{2}}{n^{\beta}\ln b_n}\Big)\leq \frac{2}{n^2},\,\,\,\,\,\,\,\,\,\, \frac{\sqrt{n}(\hat{\alpha}-1)\ln b_n}{\sqrt{b_n}} \Rightarrow N(0,\, 1)
\label{han17}
 \end{eqnarray}
for $\alpha =1$ and $b_n \leq n^{1-2\beta}/\ln n$, and
\begin{eqnarray}
   \textbf{P}\Big( |\hat{\alpha}'-2
   | \geq \frac{2\sqrt{2}}{n^{\beta}\ln b_n}\Big)\leq \frac{2}{n^2},\,\,\,\,\,\,\,\,\,\, \frac{\sqrt{n}(\hat{\alpha}'-2)\ln b_n}{b_n} \Rightarrow N(0,\, 1)
\label{han18}
 \end{eqnarray}
for $\alpha=2$ and $b^2_n \leq n^{1-\beta}/\ln n$.}
\end{theorem}

\section*{3. Numerical Simulations}
\noindent
In this section, we will compare our two estimators $\hat{\alpha} $ and $\hat{\alpha}'$ with other five estimators in the estimating error, the Type I Error and the power, including the Hill estimator[\ref{bb10}], QQ estimator[\ref{bb14}], the Moment estimator[\ref{bb6}], t-Hill estimator[\ref{bb11}, \ref{bb12}] and t-lgHill estimator[\ref{bb11}]. Since the asymptotic distribution of IPO estimator[\ref{bb13}] is unknown, we only give the estimating error of the IPO estimator in Section 3.1.
Excepting the two truncated estimators, the other six estimators can be written as
\begin{eqnarray*}
\begin{aligned}
\hat{\alpha}_H^{-1}&=&\frac{1}{m} \sum_{i=1}^{m} \log \left(\frac{X_{(n-i+1)}}{X_{(n-m)}}\right),
\end{aligned}
\end{eqnarray*}
\begin{eqnarray*}
\begin{aligned}
\hat{\alpha}_Q^{-1}&=&\frac{\sum_{j=1}^{m} \log ((m+1) / j) \log X_{(n-j+1)}-m^{-1} \sum_{j=1}^{m} \log ((m+1) / j) \sum_{j=1}^{m} \log X_{(n-j+1)}}{\sum_{j=1}^{m} \log ^{2}((m+1) / j)-m^{-1}\left(\sum_{j=1}^{m} \log ((m+1) / j)\right)^{2}},
\end{aligned}
\end{eqnarray*}
\begin{equation*}
\begin{aligned}
\hat{\alpha}_M^{-1}=M_{m, n}^{(1)}+1-\frac{1}{2}\left(1-\frac{\left(M_{m, n}^{(1)}\right)^{2}}{M_{m, n}^{(2)}}\right)^{-1},
\end{aligned}
\end{equation*}
\begin{equation*}
\begin{aligned}
\hat{\alpha}_{tH}^{-1}=\left(   \frac{1}{m} \sum_{i=1}^{m} \frac{ \mathbf{X}_{(m+1, n)}}{  \mathbf{X}_{(i, n)} } \right)^{-1}-1,
\end{aligned}
\end{equation*}
\begin{equation*}
\begin{aligned}
\hat{\alpha}_{tlH}^{-1}=\frac{M_{n}^{(2)}-\left(M_{n}^{(1)}\right)^{2}}{M_{n}^{(1)}}
\end{aligned}
\end{equation*}
and
\begin{equation*}
\hat{\alpha}_{IPO}^{-1}=-\frac{\log \left\{F_{n}^{\leftarrow}\left(\frac{3}{4}\right)+3\left[F_{n}^{\leftarrow}\left(\frac{3}{4}\right)-F_{n}^{\leftarrow}\left(\frac{1}{4}\right)\right]\right\}}{\log \hat{p}_{R}(0.25, n)},
\end{equation*}
where $X_{(1)} \leq X_{(2)} \leq \cdots \leq X_{(n)}$ are the order statistics of $X_{1} , X_{2} , \cdots , X_{n}$. $m=m(n)$, $m(n)/n \rightarrow 0$ as $n \rightarrow \infty$, $m(n)=1,2, \ldots, n-1$. $M_{m, n}^{(l)}$ is defined as
\begin{equation*}
M_{m, n}^{(l)}=\frac{1}{m} \sum_{i=1}^{m}\left(\log \frac{X_{(n-i)}}{X_{\left(n-m\right)}}\right)^{l}, \quad l=1, 2
\end{equation*}
and the definitions of $F_{n}^{\leftarrow}(p)$ and $\hat{p}_{R}(p, n)$ in detail are given in the paper[\ref{bb13}].

\subsection*{\emph{3.1. The estimating error}}
\noindent
In our next simulations, let $n=10000$ be the number of samples in each trial and set the truncated sequence be $b_{n}=n^{q}$ in the truncated estimator $\hat{\alpha}$ or  $\hat{\alpha}'$, where $q$ is the index of $b_{n}$ satisfying $0<$$\alpha$$q<1$. In order to obtain the simulation searching accuracy for the two truncated estimators, we take $k$ as the number of iterations (see  (11) and (12) in Theorem 1) such that  $\left|\hat{\alpha}-\hat{\alpha}_{k}\right| \leq \epsilon=0.001$ or $\left|\hat{\alpha}'-\hat{\alpha}_{k}'\right| \leq \epsilon=0.001$. \textbf{Remark \ref{remark2}} provides a method on how to determine the initial values $\hat{\alpha}_0$ and $\hat{\alpha}'_0$.

We first consider $\alpha \in (0, 1]$ and take the initial value $\hat{\alpha}_0=0.5$. Let $m=m(n)=\sqrt{n}=100$ for the four estimators $\hat{\alpha}_{H}$, $\hat{\alpha}_{Q}$ ,$\hat{\alpha}_{M}$ and $\hat{\alpha}_{tH}$ and set $m=m(n)={0.5*n}=5000$ for the  $\hat{\alpha}_{tlH}$.  We set $p=0.25$ for $\hat{\alpha}_{IPO}$.

The following table {\ref{tb1}} and figure {\ref{fg1}} illustrate the numerical simulation results for the \textbf{seven estimators}. All the numerical simulation results in this section were obtained using $N=10^{3}$ repetitions.
\begin{table}[H]
  \setlength{\abovecaptionskip}{0cm}
  \caption{. The estimation for different $\alpha\in (0, 1]$.}
  \setlength{\tabcolsep}{3.5mm}
    \begin{tabular*}{\hsize}{lcrclcccccccr} \arrayrulecolor{blue}\hline
     \multicolumn{3}{c}{Parameters} & \multicolumn{8}{c}{Estimation} \\\hline
   $\alpha$ &q  &k  & & $\hat{\alpha}$ &$\hat{\alpha}_{H}$ & $\hat{\alpha}_{Q}$ & $\hat{\alpha}_{M}$  & $\hat{\alpha}_{tH}$ & $\hat{\alpha}_{tlH}$ & $\hat{\alpha}_{IPO}$\\\hline
       0.10   &2.00  &11&  &0.101 &0.101 &0.097 &0.101 &0.091 &0.100 &0.157    \\
       0.20   &2.00  &6 &  &0.200 &0.201 &0.194 &0.202 &0.190 &0.201 &0.216    \\
       0.30   &2.00  &5 &  &0.301 &0.299 &0.289 &0.303 &0.287 &0.304 &0.304    \\
       0.40   &1.80  &4 &  &0.402 &0.404 &0.386 &0.408 &0.391 &0.410 &0.403    \\
       0.50   &1.70  &4 &  &0.502 &0.496 &0.481 &0.510 &0.488 &0.517 &0.503    \\
       0.60   &1.50  &5 &  &0.603 &0.595 &0.583 &0.620 &0.586 &0.624 &0.606    \\
       0.70   &1.30  &7 &  &0.704 &0.705 &0.677 &0.725 &0.686 &0.727 &0.710    \\
       0.80   &1.20  &9 &  &0.803 &0.793 &0.771 &0.831 &0.784 &0.827 &0.816    \\
       0.90   &1.10  &12&  &0.904 &0.895 &0.868 &0.943 &0.884 &0.923 &0.925    \\\hline
       &&&\textbf{AE} &\textbf{0.002} &\text{0.004} &\text{0.017} &\text{0.016}&\text{0.013} &\text{0.015}&\text{0.016}\\\hline
    \end{tabular*}
\label{tb1}
\end{table}
The \textbf{AE} in the last row of the table \ref{tb1} denotes the average of estimating error, the smaller the mean error, the better the estimator. We define $\mathbf{AE}=\sum\left|\hat{\alpha}_{Z}-\alpha\right|/9$, where $\hat{\alpha}_{Z}$ denotes one of estimators $\hat{\alpha}$, $\hat{\alpha}_{H}$, $\hat{\alpha}_{Q}$, $\hat{\alpha}_{M}$, $\hat{\alpha}_{tH}$, $\hat{\alpha}_{tlH}$ and $\hat{\alpha}_{IPO}$ for $\alpha$=0.10, 0.20, 0.30, 0.40, 0.50, 0.60, 0.70, 0.80, 0.90. 
\begin{center}
\centering
\begin{figure}[H]
    \centering
    \includegraphics[width=0.80\textwidth]{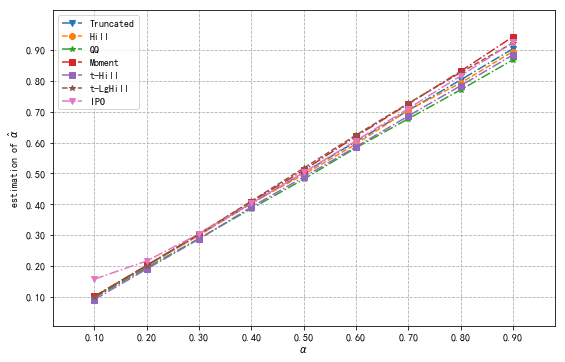}
    \caption{. The estimation for different $\alpha\in (0,1)$.}
\label{fg1}
\end{figure}
\end{center}
\vspace{-0.8cm}

It can be seen from the table \ref{tb1} and the figure \ref{fg1} that the estimating errors $\left|\hat{\alpha}-\alpha\right|$ of $\hat{\alpha}$ are smaller than that of other six estimators for $\alpha=0.20, 0.30, 0.40, 0.50, 0.60, 0.70, 0.80, 0.90$. Only for $\alpha=0.10$, the estimating error $\left|\hat{\alpha}-0.01\right|=0.001$ of $\hat{\alpha}$ is larger than that of $\hat{\alpha}_{tlH}$ since  $\left|\hat{\alpha}_{tlH}-0.10\right|=0$. When $\alpha=0.10$ or $\alpha=0.20$,  the estimating error of  $\hat{\alpha}_{IPO}$ is large than that of other six estimators. Obviously, the average of estimating error \textbf{AE} (\textbf{0.02}) of $\hat{\alpha}$ is smallest among the seven estimators.  That is, we can say that the estimator $\hat{\alpha}$ has the best performance in estimating $\alpha \in (0, 1]$ among the seven estimators.

Next, we consider $\alpha \in (1, 2]$. Take the initial  value $\hat{\alpha}'_{0}=1.5$. Similar to $\alpha \in (0, 1]$, let $n=10000$ be the number of samples. Let the truncated sequence be $b_{n}=n^{q}$ in the truncated estimator $\hat{\alpha}$, where $q$ is the index of $b_{n}$ satisfying $0<\alpha$$q<1$.
\begin{table}[H]
  \setlength{\abovecaptionskip}{0cm}
  \caption{. The estimation for different $\alpha\in (1, 2)$.}
  \setlength{\tabcolsep}{3.5mm}
    \begin{tabular*}{\hsize}{lcrclcccccccr} \arrayrulecolor{blue}\hline
     \multicolumn{3}{c}{Parameters} & \multicolumn{8}{c}{Estimation} \\\hline
 $\alpha$ &q  &k  & & $\hat{\alpha}'$ &$\hat{\alpha}_{H}$ & $\hat{\alpha}_{Q}$ & $\hat{\alpha}_{M}$  & $\hat{\alpha}_{tH}$ & $\hat{\alpha}_{tlH}$ & $\hat{\alpha}_{IPO}$\\\hline
       1.10   &0.80  &3 &  &1.108 &1.089 &1.055 &1.116 &1.075 &1.098 &0.613   \\
       1.20   &0.70  &3 &  &1.212 &1.187 &1.156 &1.283 &1.169 &1.181 &0.686   \\
       1.30   &0.65  &5 &  &1.317 &1.274 &1.245 &1.399 &1.254 &1.257 &0.766   \\
       1.40   &0.63  &8 &  &1.425 &1.365 &1.334 &1.513 &1.344 &1.327 &0.851   \\
       1.50   &0.61  &7 &  &1.537 &1.448 &1.413 &1.611 &1.426 &1.392 &0.946   \\
       1.60   &0.60  &9 &  &1.652 &1.543 &1.508 &1.736 &1.521 &1.454 &1.049   \\
       1.70   &0.60  &15&  &1.765 &1.620 &1.580 &1.846 &1.597 &1.513 &1.164   \\
       1.80   &0.58  &19&  &1.890 &1.714 &1.692 &1.985 &1.683 &1.566 &1.292   \\
       1.90   &0.45  &24&  &1.995 &1.790 &1.764 &2.166 &1.759 &1.618 &1.437   \\\hline
       &&&\textbf{AE}&\textbf{0.095} &\text{0.110} &\text{0.136} &\text{0.116}&\text{0.141} &\text{0.282}  &\text{0.463}    \\\hline
    \end{tabular*}
\label{tb2}
\end{table}
\begin{center}
\centering
\begin{figure}[H]
    \centering
    \includegraphics[width=0.80\textwidth]{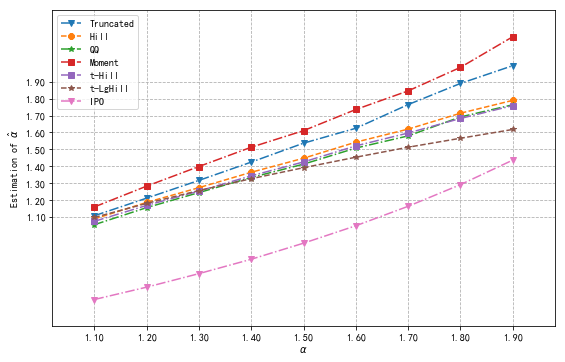}
    \caption{. The estimation for different $\alpha\in (1, 2)$.}
\label{fg2}
\end{figure}
\end{center}
\vspace{-0.8cm}

Similarly, from the table \ref{tb2} and the figure \ref{fg2} we can see that the estimating errors $\left|\hat{\alpha}'-\alpha\right|$ of $\hat{\alpha}'$ are smaller than that of other six estimators for $\alpha= 1.20, 1.30, 1.40, 1.50, 1.60, 1.70,  1.90$. Only for $\alpha=1.10$ and $\alpha=1.80$, the estimating errors $\left|\hat{\alpha}'-1.01\right|=0.008$ and $\left|\hat{\alpha}'-1.80\right|=0.090$ of $\hat{\alpha}'$ are larger than that of $\hat{\alpha}_{tlH}$ and $\hat{\alpha}_{H}$ respectively since  $\left|\hat{\alpha}_{tlH}-1.10\right|=0.002$ and $\left|\hat{\alpha}_{H}-1.80\right|=0.086$. The estimating error of  $\hat{\alpha}_{IPO}$ is large than that of other six estimators for all $\alpha$. Obviously, the average of estimating error \textbf{AE} (\textbf{0.095}) of $\hat{\alpha}'$ is smallest among the seven estimators.  That is, the estimator $\hat{\alpha}$ has the best performance in estimating $\alpha \in (1, 2]$ among the seven estimators.

In short, the two  truncated estimators $\hat{\alpha}$ and $\hat{\alpha}^{\prime}$ have the best performance in estimating $\alpha$ ($0<\alpha \leq 2$) among the seven estimators.

\begin{remark}{\label{remark3}}
The disadvantage of the two truncated estimators is that they need to know the value range of the unknown parameter $\alpha$. If we don't know whether $\alpha$ is included in interval $(0, 1]$ or interval $(1, 2]$, we may take the initial values $\hat{\alpha}_{0}$ and $\hat{\alpha}'_{0}$ according to the method in  \textbf{Remark \ref{remark2}}
\end{remark}

\subsection*{\emph{3.2. The rejection regions and the Type I Error}}
\noindent
In order to get the Type I Error, we consider the rejected regions of these estimators except the IPO estimator since we do not know the asymptotic distribution of $ \hat{\alpha}_{IPO}$. Let $H_{0}$ and $H_{1}$ denote the original hypothesis and the alternative hypothesis, respectively, that is,
\begin{eqnarray*}
\text{Original hypothesis} \,\, H_0: \,\, \alpha=\alpha_0,\,\,\,\,\,\,\,\,\,\,\,             \text{ alternative hypothesis}\,\,  H_1:   \,\, \alpha\neq \alpha_0,
\end{eqnarray*}
where $0<\alpha_{0}< 1$ or $1<\alpha_{0}< 2$. Let the confidence level be $0.95$ and $\beta=0$ in the \textbf{Theorem \ref{theorem3}}. By using the inequalities (\ref{han15}) and (\ref{han16}) of the \textbf{Theorem \ref{theorem3}}, we have
\begin{eqnarray*}
\textbf{P}\Big( \frac{\sqrt{n}\sqrt{\alpha_0(2-\alpha_0)}|\hat{\alpha}-\alpha_0|\ln b_n}{(1-\alpha_0)b_n^{\alpha_0/2}} \leq 1.96\Big)\approx 2\Phi(1.96)-1=0.95
\end{eqnarray*}
for $0<\alpha_0<1$ and
\begin{eqnarray*}
\textbf{P}\Big( \frac{\sqrt{n}\sqrt{\alpha_0(4-\alpha_0)}|\hat{\alpha}'-\alpha_0|\ln b_n}{(2-\alpha_0)b_n^{\alpha_0/2}} \leq 1.96\Big)\approx 2\Phi(1.96)-1=0.95
\end{eqnarray*}
for $1<\alpha_0<2$. Therefore, we can get two rejection regions $R_T$ and $R_T'$ in the following
\begin{eqnarray*}
R_T=\{x:  \frac{\sqrt{n}\sqrt{\alpha_0(2-\alpha_0)}|x-\alpha_0|\ln b_n}{(1-\alpha_0)b_n^{\alpha_0/2}} > 1.96\}
\end{eqnarray*}
for $0<\alpha_0<1$ and
\begin{eqnarray*}
R_T=\{x:  \frac{\sqrt{n}\sqrt{\alpha_0(4-\alpha_0)}|x-\alpha_0|\ln b_n}{(2-\alpha_0)b_n^{\alpha_0/2}} > 1.96\}
\end{eqnarray*}
for $1<\alpha_0<2$.

Since the five estimators, $\hat{\alpha}_{H}$, $\hat{\alpha}_{Q}$, $\hat{\alpha}_{M}$, $\hat{\alpha}_{tH}$ and  $\hat{\alpha}_{tlH}$ satisfy
\begin{equation*}
\alpha_{0} \sqrt{m}\left(\hat{\alpha}_{H}^{-1}-\alpha_{0}^{-1}\right) \stackrel{d}{\longrightarrow} N(0,1),
\end{equation*}
\begin{equation*}
\alpha_{0} \sqrt{m / 2}\left(\hat{\alpha}_{Q}^{-1}-\alpha_{0}^{-1}\right) \stackrel{d}{\longrightarrow} N(0,1),
\end{equation*}
\begin{equation*}
\frac{\alpha_{0} \sqrt{m}}{\sqrt{1+\alpha_{0}^{2}}}\left(\hat{\alpha}_{M}^{-1}-\alpha_{0}^{-1}\right) \stackrel{d}{\longrightarrow} N(0,1),
\end{equation*}
\begin{equation*}
\frac{\alpha_{0}\sqrt{\alpha_{0}(\alpha_{0}+2)}{\sqrt{m}}}{1+\alpha_{0}}\left( \hat{\alpha}_{tH}^{-1}-\alpha_{0}^{-1} \right) \stackrel{d}{\longrightarrow} N(0,1),
\end{equation*}
and
\begin{equation*}
\frac{\alpha_{0}{\sqrt{m}}}{2\sqrt{2}}\left( \hat{\alpha}_{tlH}^{-1}-\alpha_{0}^{-1} \right) \stackrel{d}{\longrightarrow} N(0,1),
\end{equation*}
we can similarly  get the five rejection regions $R_{H}$, $R_{Q}$, $R_{M}$,  $R_{tH}$ and $R_{tlH}$ in the following with the confidence level $0.95$ respectively
\begin{equation*}
R_{H} =\left\{x: \alpha_{0} \sqrt{m}\left|x^{-1}-\alpha_{0}^{-1}\right|>1.96\right\},
\end{equation*}
\begin{equation*}
R_{Q} =\left\{x: \alpha_{0} \sqrt{m / 2}\left|x^{-1}-\alpha_{0}^{-1}\right|>1.96\right\},
\end{equation*}
\begin{equation*}
R_{M}=\left\{x: \frac{\alpha_{0} \sqrt{m}}{\sqrt{1+\alpha_{0}^{2}}}\left|x^{-1}-\alpha_{0}^{-1}\right|>1.96\right\},
\end{equation*}
\begin{equation*}
R_{tH}=\left\{x: \frac{\alpha_{0}\sqrt{\alpha_{0}(\alpha_{0}+2)}{\sqrt{m}}}{1+\alpha_{0}}\left|x^{-1}-\alpha_{0}^{-1}\right|>1.96\right\},
\end{equation*}
and
\begin{equation*}
R_{tlH}=\left\{x: \frac{\alpha_{0}{\sqrt{m}}}{2\sqrt{2}}\left|x^{-1}-\alpha_{0}^{-1}\right|>1.96\right\}.
\end{equation*}

Similar to the Section 3.1,  we first consider $\alpha_{0} \in (0, 1)$ and set the initial  value $\hat{\alpha}_{0}=0.5$. Let $m=m(n)=\sqrt{n}=100$ for the four estimators $\hat{\alpha}_{H}$, $\hat{\alpha}_{Q}$, $\hat{\alpha}_{M}$ and $\hat{\alpha}_{tH}$ and set $m=m(n)={0.5*n}=5000$ for  $\hat{\alpha}_{tlH}$.
\begin{table}[H]
  \setlength{\abovecaptionskip}{0cm}
  \caption{. The Type I Error for different $\alpha_{0} \in (0,1)$.}
  \setlength{\tabcolsep}{4.5mm}
    \begin{tabular*}{\hsize}{lcrclccccccr} \arrayrulecolor{blue}\hline
     \multicolumn{3}{c}{Parameters} & \multicolumn{7}{c}{Type I Error} \\\hline
 $\alpha_0$ &q  &k  & & $\hat{\alpha}$  &$\hat{\alpha}_{H}$ & $\hat{\alpha}_{Q}$ & $\hat{\alpha}_{M}$  & $\hat{\alpha}_{tH}$ & $\hat{\alpha}_{tlH}$ \\\hline
       0.10   &2.00  &11&  &0.038 &0.058 &0.074 &0.053 &0.171 &0.009  \\
       0.20   &2.00  &6 &  &0.052 &0.061 &0.073 &0.050 &0.106 &0.006  \\
       0.30   &2.00  &5 &  &0.046 &0.062 &0.071 &0.052 &0.107 &0.017  \\
       0.40   &1.80  &4 &  &0.050 &0.055 &0.065 &0.057 &0.069 &0.041  \\
       0.50   &1.70  &4 &  &0.055 &0.050 &0.074 &0.042 &0.076 &0.068  \\
       0.60   &1.50  &5 &  &0.055 &0.061 &0.068 &0.059 &0.088 &0.095  \\
       0.70   &1.30  &7 &  &0.046 &0.050 &0.067 &0.054 &0.069 &0.053  \\
       0.80   &1.20  &9 &  &0.026 &0.053 &0.072 &0.047 &0.071 &0.051  \\
       0.90   &1.10  &12&  &0.098 &0.053 &0.082 &0.070 &0.074 &0.020  \\\hline
       &&&\textbf{AT} &\textbf{0.052} &\text{0.056} &\text{0.072} & \text{0.054}& \text{0.092}& \text{0.040} \\\hline
    \end{tabular*}
\label{tb3}
\end{table}

Like the definition of the average of estimating error $\textbf{AE}$ we can similarly define the average of Type I Error, $\mathbf{AT}=\sum{Type I Error}/9$ under the confidence level $\mathbf{0.95}$.  The closer the average  ($\mathbf{0.050}$)  of the Type I Error, the better the estimator.

\begin{center}
\centering
\begin{figure}[H]
    \centering
    \includegraphics[width=0.80\textwidth]{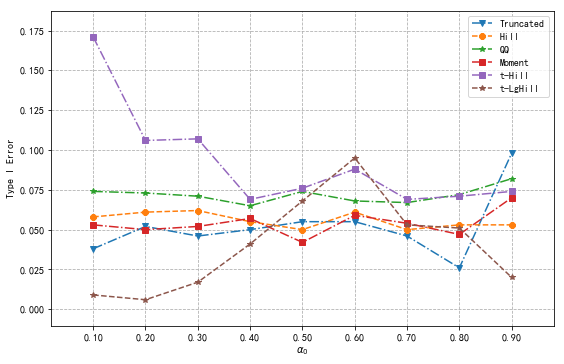}
 \caption{. The Type I Error for different $\alpha_{0} \in (0,1).$}
\label{fg3}
\end{figure}
\end{center}
\vspace{-0.8cm}

From the table \ref{tb3} and the figure \ref{fg3}  we can see that the value $ \textbf{AT} $ of $\hat{\alpha}$ is closer to $0.050$ than that of the other five estimators. Thus, it could be said that the truncated estimator $\hat{\alpha}$ is better than other five estimators for estimating $\alpha_{0} \in (0, 1)$.

Next, we consider $\alpha_{0} \in (1, 2)$ and set the initial  value $\hat{\alpha}'_{0}=1.5$. Similar to $\alpha_{0} \in(0, 1)$, let $n=10000$ be the number of samples and the truncated sequence be $b_{n}=n^{q}$ in the truncated estimator $\hat{\alpha}'$, where $q$ is the index of $b_{n}$ satisfying $0<\alpha_{0} q<1$.
\begin{table}[H]
  \setlength{\abovecaptionskip}{0cm}
  \caption{. The Type I Error for different $\alpha_{0} \in (1,2)$.}
  \setlength{\tabcolsep}{4.5mm}
    \begin{tabular*}{\hsize}{lcrclccccccr} \arrayrulecolor{blue}\hline
 \multicolumn{3}{c}{Parameters} &\multicolumn{7}{c}{Type I Error}  \\\hline
	   $\alpha_0$ &q  &k &  & $\hat{\alpha}'$ &$\hat{\alpha}_{H}$ & $\hat{\alpha}_{Q}$ & $\hat{\alpha}_{M}$ & $\hat{\alpha}_{TH}$ & $\hat{\alpha}_{tlH}$ \\\hline
       1.10   &0.80  &3&  &0.043 &0.063 &0.089 &0.060 &0.078 &0.004   \\
       1.20   &0.70  &3&  &0.041 &0.052 &0.064 &0.068 &0.072 &0.007   \\
       1.30   &0.65  &5&  &0.049 &0.054 &0.071 &0.057 &0.090 &0.048   \\
       1.40   &0.63  &8&  &0.053 &0.079 &0.081 &0.068 &0.095 &0.190   \\
       1.50   &0.61  &7&  &0.067 &0.101 &0.101 &0.049 &0.119 &0.489   \\
       1.60   &0.60  &9&  &0.090 &0.079 &0.090 &0.062 &0.094 &0.803   \\
       1.70   &0.60  &15& &0.094 &0.104 &0.096 &0.056 &0.135 &0.964   \\
       1.80   &0.58  &19& &0.093 &0.097 &0.084 &0.059 &0.128 &0.999   \\
       1.90   &0.45  &24& &0.199 &0.125 &0.112 &0.069 &0.151 &1.000   \\\hline
       &&&\textbf{AT} &\text{0.081} &\text{0.084} &\text{0.088} &\textbf{0.061} &\text{0.107}&\text{0.500}    \\\hline
    \end{tabular*}
\label{tb4}
\end{table}
\begin{center}
\centering
\begin{figure}[H]
    \centering
    \includegraphics[width=0.80\textwidth]{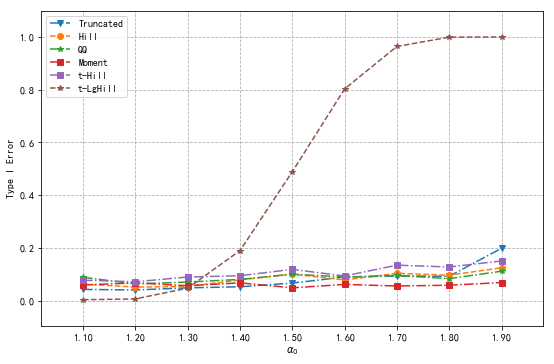}
\caption{. The Type I Error for different $\alpha_{0}\in(1,2)$.}
\label{fg4}
\end{figure}
\end{center}
\vspace{-0.8cm}

From the table \ref{tb4} and the figure \ref{fg4} we can see that the value $\mathbf{A T}=0.081$ of $\hat{\alpha}'$ is closer to $0.050$ than that of other four estimators except the Moment estimator $\hat{\alpha}_{M}$ since the average Type I Error of $\hat{\alpha}_{M}$ is $\textbf{0.061}$.

\subsection*{\emph{3.3. Power of estimator}}
\noindent
In this section we consider the power of estimator, that is, the probability of correctly rejecting the original hypothesis under the confidence level $0.95$. Consider two original hypothesises $\alpha_{0}=0.6$ and $\alpha_{0}=1.40$ respectively. Take $b_{n}=n^{1.5}$, $n=10000$ and consider several different tail indices $\alpha^*=0.60, 0.64, 0.68, 0.72, 0.76, 0.80, 0.84, 0.88, 0.92$, we can get the corresponding estimators  $\hat{\alpha}^{*}$, $\hat{\alpha}_{H}^{*}, \hat{\alpha}_{Q}^{*}, \hat{\alpha}_{M}^{*}$, $\hat{\alpha}_{tH}^{*}$ and  $\hat{\alpha}_{tlH}^{*}$. We can similarly define the average power as $\mathbf{A P}=\sum{\hat{\alpha}_{P}^{*}} / 9$, where $\hat{\alpha}_{P}^{*}$ denotes the power of  $\hat{\alpha}^{*}$, $\hat{\alpha}_{H}^{*}, \hat{\alpha}_{Q}^{*}, \hat{\alpha}_{M}^{*}$, $\hat{\alpha}_{tH}^{*}$ and  $\hat{\alpha}_{tlH}^{*}$.

\begin{table}[H]
  \setlength{\abovecaptionskip}{0cm}
  \caption{. The power for $\alpha_{0}=0.60$ with different $\alpha^{*}$.}
  \setlength{\tabcolsep}{5.5mm}
    \begin{tabular*}{\hsize}{lrclccccccr} \arrayrulecolor{blue}\hline
 \multicolumn{2}{c}{Parameters} &\multicolumn{7}{c}{Power}  \\\hline
        $\alpha^*$ &k&  &$\hat{\alpha}^{*}$ &$\hat{\alpha}_{H}^{*}$ &$\hat{\alpha}_{Q}^{*}$ &$\hat{\alpha}_{M}^{*}$ &$\hat{\alpha}_{tH}^{*}$  &$\hat{\alpha}_{tlH}^{*}$\\\hline
        $0.60$ &5  &  &0.055 &0.056 &0.062 &0.049 &0.064 &0.076  \\
        $0.64$ &5  &  &0.137 &0.067 &0.029 &0.086 &0.049 &0.777  \\
        $0.68$ &6  &  &0.552 &0.194 &0.037 &0.215 &0.083 &0.999  \\
        $0.72$ &7  &  &0.841 &0.330 &0.108 &0.342 &0.147 &1.000  \\
        $0.76$ &8  &  &0.967 &0.536 &0.179 &0.509 &0.269 &1.000  \\
        $0.80$ &5  &  &0.990 &0.760 &0.295 &0.657 &0.440 &1.000  \\
        $0.84$ &9  &  &0.995 &0.875 &0.407 &0.772 &0.604 &1.000  \\
        $0.88$ &7  &  &1.000 &0.945 &0.558 &0.893 &0.709 &1.000  \\
        $0.92$ &9  &  &1.000 &0.985 &0.628 &0.925 &0.835 &1.000  \\\hline
        &&\textbf{AP}&\text{0.726} &\text{0.528} &\text{0.256} &\text{0.494} &\text{0.356} &\textbf{0.872} \\\hline
    \end{tabular*}
\label{tb5}
\end{table}

The  table \ref{tb5}  above and the following figure \ref{fg5}  illustrate  the power and average power $\textbf{AP}$ of six estimators, $\hat{\alpha}^{*}$, $\hat{\alpha}_{H}^{*}$, $\hat{\alpha}_{Q}^{*}$, $\hat{\alpha}_{M}^{*}$, $\hat{\alpha}_{tH}^{*}$ and $\hat{\alpha}_{tlH}^{*}$. It can be seen that the average power $\textbf{AP}$ of the truncated estimator $\hat{\alpha}^{*}$ is 0.726 which is larger than that of other four estimators except the t-lgHill estimator since  the average power of $\hat{\alpha}_{tlH}^{*}$ is $\textbf{0.872}$.

\begin{center}
\centering
\begin{figure}[H]
  \centering
  \includegraphics[width=0.80\textwidth]{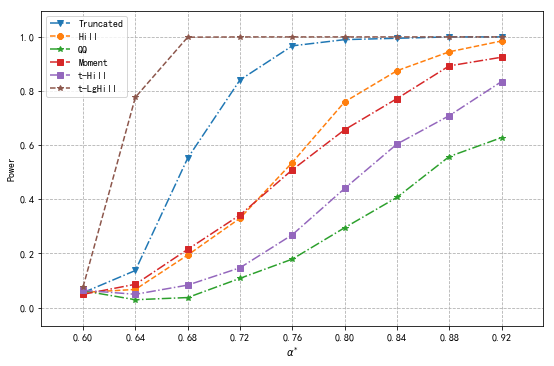}
  \caption{. The power for $\alpha_{0}=0.60$ with different $\alpha^{*}$.}
\label{fg5}
\end{figure}
\end{center}
\vspace{-0.8cm}

Next we consider $\alpha_{0}=1.40$. It can be seen from the following table \ref{tb6} and figure \ref{fg6} that the power of $\hat{\alpha}^{* \prime}$ is larger than that of other five estimators respectively for $\alpha^{* \prime}=1.48, 1.56, 1.64, 1.72$, $ 1.80, 1.88, 1.96$, and the average power $\textbf{AP}$ ($\textbf{0.781}$) of $\hat{\alpha}^{* \prime}$ is the largest among all six estimators.

\begin{table}[H]
  \setlength{\abovecaptionskip}{0cm}
  \caption{. The power for $\alpha_{0}=1.40$ with different $\alpha^{* \prime}$.}
  \setlength{\tabcolsep}{5.5mm}
    \begin{tabular*}{\hsize}{lrclccccccr} \arrayrulecolor{blue}\hline
	$\alpha^{* \prime}$ &k &  &$\hat{\alpha}^{* \prime}$ &$\hat{\alpha}_{H}^{* \prime}$ &$\hat{\alpha}_{Q}^{* \prime}$ &$\hat{\alpha}_{M}^{* \prime}$ &$\hat{\alpha}_{tH}^{* \prime}$  &$\hat{\alpha}_{tlH}^{* \prime}$\\\hline
        $1.40$ &8   &  &0.068 &0.066 &0.084 &0.067 &0.096 &0.170  \\
        $1.48$ &8   &  &0.383 &0.038 &0.040 &0.088 &0.044 &0.006  \\
        $1.56$ &7   &  &0.820 &0.076 &0.036 &0.120 &0.055 &0.008  \\
        $1.64$ &8   &  &0.978 &0.164 &0.058 &0.181 &0.107 &0.125  \\
        $1.72$ &11  &  &0.999 &0.256 &0.080 &0.207 &0.176 &0.578  \\
        $1.80$ &12  &  &1.000 &0.414 &0.142 &0.331 &0.264 &0.925  \\
        $1.88$ &18  &  &1.000 &0.562 &0.196 &0.375 &0.389 &0.995  \\
        $1.96$ &26  &  &1.000 &0.678 &0.254 &0.423 &0.518 &1.000  \\\hline
        &&\textbf{AP}&\textbf{0.781} &\text{0.282} &\text{0.111} &\text{0.224} &\text{0.206} &\text{0.476} \\\hline
    \end{tabular*}
\label{tb6}
\end{table}

\begin{center}
\centering
\begin{figure}[H]
  \centering
  \includegraphics[width=0.80\textwidth]{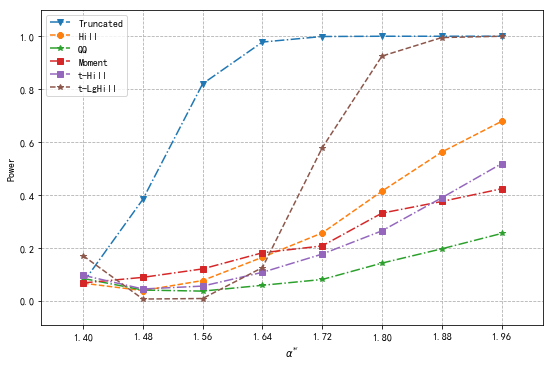}
  \caption{. The Power of $\alpha_{0}^{\prime}=1.40$ for different $\alpha^{* \prime}$.}
\label{fg6}
\end{figure}
\end{center}
\vspace{-0.8cm}

The second largest average powers of $\hat{\alpha}^{*}$ and the largest average powers of $\hat{\alpha}^{* \prime}$ respectively in table \ref{tb5} and table \ref{tb6} mean that the two truncated estimators have a more robust performance than that of other five estimators on the whole.

\section*{4. Conclusion}

In order to make up for the shortcomings of existing estimation methods, we present a new method of truncated estimation to estimate the tail index of the extremely heavy-tailed distributions with infinite mean or variance. By using the truncated sample mean $\hat{\mu}_{n}$, the truncated sample second order moment $\hat{\nu}_{n}^{2}$ and the two recursive estimators in equation (\ref{han9}) and (\ref{han10}), we can obtain the two truncated estimators $\hat{\alpha}$ and $\hat{\alpha}^{\prime}$ respectively for $\alpha\in(0, 1]$ and $\alpha\in(1, 2]$. We not only give the rate of strong consistency convergence of the two truncated estimators, but also prove that their asymptotic distributions are normal. Moreover, among  all six estimators, the numerical simulation results show that the two truncated estimators have the smallest average estimating error, the truncated estimator $\hat{\alpha}$  has the closest  average (0.05) of  Type I Error and the truncated estimator $\hat{\alpha}'$ has the largest average power. In short,  the performance of the two new truncated estimators is quite good on the whole.

\section*{Acknowledgments}
\noindent
The authors are grateful to the referees for their careful reading of this paper and valuable comments.

\section*{Declaration of interest statement}
\noindent
The authors declare that they have no conflict of interest.

\section*{Availability of data and materials}
\noindent
The datasets used and/or analysed during the current study are available from the corresponding author on reasonable request.


~~~\\
~~~\\

\noindent
{\Large{\textbf{Appendix: Proofs Theorems}}}\\
~~~\\
~~~\\

\noindent
\textbf{Proof of Theorem 1. } Let $H_1(x)=x-G_1(x)$ for $0<x\leq 1$ and $H_2(y)=y-G_2(y)$ for $0<y\leq 2$. Let
\begin{eqnarray*}
&&0=H'_1(x)=1-G'_1(x)=1-\frac{\hat{\mu}(b_n)}{ \ln b_n[x(1-x)\hat{\mu}(b_n)+x^2]},\,\,\,\, 0<x \leq 1, \\
&&0=H'_2(y)=1-G'_2(y)=1-\frac{2\hat{\nu}^2(b_n)}{\ln b_n[y(2-y)\hat{\nu}^2(b_n)+y^2]}, \,\,\,\, 0<y\leq 2.
\end{eqnarray*}
It follows that both  $H'_1(x)=0$ and $H'_2(y)=0$ have two  real roots, $x_{1, 2}$ and $y_{1, 2}$, respectively, i.e.
\begin{eqnarray*}
x_{1, 2}=\frac{1\pm \sqrt{1-4(1-1/\hat{\mu}(b_n))/\ln b_n}}{2(1-1/\hat{\mu}(b_n))},\,\,\,\,y_{1, 2}=\frac{1\pm \sqrt{1-2(1-1/\hat{\nu}^2(b_n))/\ln b_n}}{(1-1/\hat{\nu}^2(b_n))}.
\end{eqnarray*}
for large $n$ such that $ \ln b_n \geq 4$, $\hat{\mu}(b_n)>1$ and $\hat{\nu}^2(b_n)>1$.  Since
\begin{eqnarray*}
G''_1(x)=-\frac{\hat{\mu}(b_n)[2x+(1-2x)\hat{\mu}(b_n)]}{\log b_n [x(1-x)\hat{\mu}(b_n)+x^2]^2},
\end{eqnarray*}
it follows that $G''_1(x)<0$ for $0<x\leq 1/2$ and $G''_1(x)>0$ for $\hat{\mu}(b_n)> 2x/(2x-1)$ and $1/2\,<x<1$. Hence,  $H_1(x)=x-G_1(x)$ is monotonically increasing for $x_1<x<x_2$ since $G'_1(x)>0$ for $0<x<1$.

Let $\mu_n(x_1)=x_1(1-x_1)^{-1}(b_n^{1-x_1}-1)$. Note that $x_1=(1+o(1))\ln b_n(1+\ln^{-1} b_n)$ for large $n$. By using $1=G'_1(x_1)$, (\ref{han5}) and the probability of $(\hat{\mu}(b_n)-\mu_n(x_1))/\mu_n(x_1)$ in (\ref{han13}) of the theorem 2, we can get that
\begin{eqnarray*}
H_1(x_1)&=&x_1-1+\frac{\ln [\hat{\mu}(b_n)/\ln b_n]-2\ln x_1}{\ln b_n} \\
&=&x_1-1+\frac{\ln [\mu_n(x_1)/\ln b_n]-2\ln x_1+ \ln [(\hat{\mu}(b_n)-\mu_n(x_1))/\mu_n(x_1)+1]}{\ln b_n}\\
&\leq &\frac{-\ln (1-x_1)-\ln \ln b_n-\ln x_1-2/\sqrt{\ln b_n}}{\ln b_n}\\
&= & (1+o(1))\frac{2(x_1-\sqrt{x_1})}{\ln b_n}<0
\end{eqnarray*}
with high probability $1-n^{-2}$ for large $n$.  On the other hand, $H_1'(x_2)=0$, $H_1(1)=1-G_1(1)=0$ and $H_1'(1)= 1-G'_1(1)<0$ since $\hat{\mu}(b_n)(\ln b_n)^{-1}>1$ for  large $n$, it follows that $H_1(x_2)>0$ for large $n$. Thus, $H_1(x)=0$ has an unique root $\hat{\alpha}\in (x_1,\, x_2)$ such that $H(\hat{\alpha})=0$, i.e. $\hat{\alpha}=G_1(\hat{\alpha})$ for large $n$. Note that $(x_1,\,\ x_2)\to (0,\,\,1)$ as $n \to \infty$, therefore, $H_1(x)=0$ has an unique root $\hat{\alpha}\in (0,\, 1)$ for large $n$.

Similarly, from
\begin{eqnarray*}
G''_2(y)=-\frac{4\hat{\nu}^2(b_n)[y+(1-y)\hat{\nu}^2(b_n)]}{\log b_n [y(2-y)\hat{\nu}^2(b_n)+y^2]^2},
\end{eqnarray*}
it follows that $G''_2(y)<0$ for $0<y\leq 1$ and $G''_2(y)>0$ for $\hat{\nu}^2(b_n)> y/(y-1)$ and $1<y<2$. Hence,  $H_2(y)=y-G_2(y)$ is monotonically increasing for $y_1<y<y_2$ since $G'_2(y)>0$ for $0<y<2$.

Let $\varepsilon_n = \ln^{-1} b_n$. By equation(\ref{han6}) and the probability of $(\hat{\nu}^2(b_n)-\nu^2_n)/\nu^2_n$ in (\ref{han14}) of the theorem 2,  we can get that
\begin{eqnarray*}
H_2(1+\varepsilon_n)&=&\nu^{-2}_n-1+\frac{\ln[\frac{1-\varepsilon_n}{1+\varepsilon_n}\hat{\nu}^2_n +1]}{\ln b_n} \\
&=&\nu^{-2}_n-1+\frac{\ln \nu^2_n +\ln [\frac{1-\varepsilon_n}{1+\varepsilon_n}]+\ln [1+\frac{\hat{\nu}^2(b_n)-\nu^2_n}{\nu^2_n}+\frac{1+\varepsilon_n}{\nu^2_n(1-\varepsilon_n)}]}{\ln b_n}\\
&\leq &\frac{-\frac{2}{\sqrt{\ln b_n}}+\frac{1+\varepsilon_n}{\nu^2_n(1-\varepsilon_n)}}{\ln b_n} <0
\end{eqnarray*}
with high probability $1-n^{-2}$ for large $n$ since $\nu^2_n = (1+\varepsilon_n)(1-\varepsilon_n)^{-1}(b_n^{1-\varepsilon_n}-1)$.  On the other hand, since $H_2'(y_2)=0$, $H_2(2)=0$ and $H_2'(2)<0$, it follows that $H_2(y_2)>0$. Thus, $H_2(y)=0$ has an unique root $\hat{\alpha}'\in (1+\varepsilon_n,\, y_2)$ such that $H(\hat{\alpha}')=0$, i.e. $\hat{\alpha}'=G_2(\hat{\alpha}')$ for large $n$. Note that $(1+\varepsilon_n,\,\ y_2)\to (1,\,\,2)$ as $n \to \infty$, therefore, $H_2(y)=0$ has an unique root $\hat{\alpha}\in (1,\, 2)$ for large $n$.

Note that the functions $H_1(x)$ ($0<x<1$) is monotonically increasing for large $n$.  Let $0<\hat{\alpha}_0<\hat{\alpha}<1$. Since $H_1(\hat{\alpha})=0$, it follows $H_1(\hat{\alpha}_0)< H_1(\hat{\alpha})=0$ and therefore, $\hat{\alpha}_0< G_1(\hat{\alpha}_0)=\hat{\alpha}_1$. Through step-by-step iteration, we can get $\hat{\alpha}_k\nearrow \hat{\alpha}$. Furthermore, by (\ref{han3}), (\ref{han5}) and (\ref{han9}) we have
\begin{eqnarray*}
     \hat{\alpha}-\hat{\alpha}_k &=&\frac{1}{\ln b_n}\ln\Big(1+\frac{\hat{\mu}(b_n)(\hat{\alpha}-\hat{\alpha}_{k-1})}{\hat{\alpha}_{k-1}(\hat{\alpha}+\hat{\mu}(b_n))(1-\hat{\alpha})}\Big)\\
          &\leq & \frac{\hat{\mu}(b_n)(\hat{\alpha}-\hat{\alpha}_{k-1})}{\hat{\alpha}_{k-1}[\hat{\alpha}+\hat{\mu}(b_n)(1-\hat{\alpha})]\ln b_n}\\
          &\leq &\frac{\hat{\alpha}-\hat{\alpha}_{k-1}}{\hat{\alpha}_{k-1}(1-\hat{\alpha})\ln b_n}\\
          &\leq &\frac{\hat{\alpha}-\hat{\alpha}_{0}}{[\hat{\alpha}_{0}(1-\hat{\alpha})\ln b_n]^k}
\end{eqnarray*}
 since $\ln(1+x)|\leq x$ for $x\geq 0$. For $0<\hat{\alpha}< \hat{\alpha}_0<1$, we can similarly get that  $\hat{\alpha}_k \searrow \hat{\alpha}$ and
\begin{eqnarray*}
     \hat{\alpha}_k-\hat{\alpha} \leq \frac{\hat{\alpha}_{0}-\hat{\alpha}}{[\hat{\alpha}(1-\hat{\alpha})\ln b_n]^k}.
 \end{eqnarray*}
Hence, (\ref{han11}) holds.

 By the same method above we can get that  $\hat{\alpha}'_k\nearrow \hat{\alpha}'$ for $1<\hat{\alpha}'_0< \hat{\alpha}'<2$,  $\hat{\alpha}'_k \searrow \hat{\alpha}'$ for $1< \hat{\alpha}'<\hat{\alpha}'_0<2$  and the inequality (\ref{han12}) holds since $H_2(y)$ ($1<y<2$) is monotonically increasing for large $n$ and
 \begin{eqnarray*}
 \hat{\alpha}'-\hat{\alpha}'_k &=&\frac{1}{\ln b_n}\ln\Big(1+\frac{2\hat{\nu}^2(b_n)(\hat{\alpha}'-\hat{\alpha}'_{k-1})}{\hat{\alpha}'_{k-1}[\hat{\alpha}'+\hat{\nu}^2(b_n)(2-\hat{\alpha}')]}\Big),\,\,\,\,\hat{\alpha}'_0< \hat{\alpha}'\\
 \hat{\alpha}'_k-\hat{\alpha}' &=&\frac{1}{\ln b_n}\ln\Big(1+\frac{2\hat{\nu}^2(b_n)(\hat{\alpha}'_{k-1}-\hat{\alpha}')}{\hat{\alpha}'[\hat{\alpha}'_{k-1}+\hat{\nu}^2(b_n)(2-\hat{\alpha}'_{k-1})]}\Big),\,\,\,\, \hat{\alpha}'<\hat{\alpha}'_0.
 \end{eqnarray*}
 This completes the proof.\\

\noindent
 \textbf{Proof of Theorem 2. } Let $\varepsilon= 2\mu_n/(n^{\beta}\sqrt{\ln b_n})$. Note that $\textbf{P}(|X_k(b_n)-\mu_n|>b_n)=0,$ $b^{\alpha}_n \ln b_n\leq \alpha n^{1-2\beta}/\ln n$, $2-\alpha\geq 2(1-\alpha)^2$, $\mu_n=\alpha(1-\alpha)^{-1}(b_n^{1-\alpha}-1)$ and $Var(X_1(b_n))\leq \nu^2_n=\alpha(2-\alpha)^{-1}(b_n^{2-\alpha}-1)$. By the Bernstein inequality, we can get that
\begin{eqnarray*}
     \textbf{P}\Big(|\sum_{k=1}^n(X_k(b_n)-\mu_n)|\geq n\varepsilon\Big)&\leq & 2\exp\{-\frac{n\varepsilon^2}{2Var(X_1(b_n))+2b_n\varepsilon/3}\}\\
     &\leq & 2\exp\{-\frac{2n^{1-2\beta} \mu^2_n}{\nu^2_n \ln b_n  +2b_n\mu_n\sqrt{\ln b_n}/3n^{\beta}}\}\\
     &\leq &2\exp\{-\frac{ \alpha n^{1-2\beta}  (2-\alpha)b_n^{2-2\alpha}}{(1-\alpha)^2b_n^{2-\alpha}\ln b_n}\}\\
     &\leq &2\exp\{-\frac{2 \alpha n^{1-2\beta}}{b_n^{\alpha} \ln b_n }\}\leq \frac{2}{n^{2}}
 \end{eqnarray*}
for $0<\alpha <1$ and large $n$. The inequality above holds also for $\alpha=1$ since $\mu_n=\ln b_n$ and $\nu^2_n=b_n -1$ in this case. Hence, the inequality in (\ref{han13}) holds for $0<\alpha \leq 1$.
Note that
\begin{eqnarray*}
    \frac{\sum_{k=1}^n\textbf{E}(|X_k(b_n)-\mu_n|^3)}{[nVar(X_1(b_n))]^{3/2}}=O\Big(\frac{b_n^{\alpha}}{n}\Big)^{1/2} \rightarrow 0.
\end{eqnarray*}
It follows from the Lyapunov central limit theorem that
\begin{equation*}
\frac{\sqrt{(2-\alpha)} \sqrt{n}\left(\hat{\mu}_{n}-\mu_{n}\right)}{\alpha b_{n}^{1-\alpha / 2}}=(1+o(1)) \frac{n\left(\hat{\mu}_{n}-\mu_{n}\right)}{\sqrt{n \operatorname{Var}\left(X_{1}\left(b_{n}\right)\right)}} \Rightarrow N(0,1)
\end{equation*}
for $0<\alpha \leq 1$ since $Var(X_1(b_n))=(1+o(1))\alpha (2-\alpha)^{-1}b_n^{2-\alpha}$ for large $n$, i.e., (13) holds.

By the same method we can get (\ref{han14}) for $1<\alpha \leq 2$ since  $ Var(X^2_1(b_n))=(1+o(1))\alpha (4-\alpha)^{-1}b_n^{4-\alpha}$ for large $n$. It completes the proof.\\

\noindent
\textbf{Proof of Theorem 3.}  It follows from (\ref{han3}) and (\ref{han5}) that
\begin{eqnarray*}
    \alpha- \hat{\alpha}&=&\frac{1}{\ln b_n}\ln \Big(1+(1+o(1))\frac{\alpha -\hat{\alpha}}{\hat{\alpha}(1-\alpha)}+(1+o(1))\frac{\alpha(1-\hat{\alpha})}{\hat{\alpha}(1-\alpha)}\frac{(\hat{\mu}_n-\mu_n)}{\mu_n}\Big)\\
          &=&(1+o(1))\frac{1}{\ln b_n}\frac{(\hat{\mu}_n-\mu_n)}{\mu_n}
\end{eqnarray*}
for $0<\alpha <1$ and large $n$. Hence, by (\ref{han13}) of the theorem 2, we have
\begin{eqnarray*}
   \textbf{P}\Big( |\hat{\alpha}-\alpha| \geq \frac{2}{n^{\beta}\ln b_n\sqrt{\ln b_n}}\Big)= \textbf{P}\Big(|\frac{\hat{\mu}_n -\mu_n}{\mu_n}|\geq \frac{2(1+o(1))}{n^{\beta}\sqrt{\ln b_n}}\Big)\leq \frac{2}{n^2}
 \end{eqnarray*}
for $0<\alpha <1$ and large $n$. Furthermore, by $(\alpha- \hat{\alpha})\mu_n \ln b_n =(1+o(1))(\hat{\mu}_n-\mu_n)$ and (\ref{han13}), we can get that
\begin{eqnarray*}
      \frac{\sqrt{\alpha(2-\alpha)}}{1-\alpha}\frac{\sqrt{n}(\hat{\alpha}-\alpha)\ln b_n}{b_n^{\alpha/2}} (1+o(1)) &=&  \frac{\sqrt{n} (\alpha- \hat{\alpha})\mu_n \ln b_n} {\sqrt{Var(X_1(b_n))}}\\
      &=&\frac{\sqrt{n}(1+o(1))(\hat{\mu}_n-\mu_n)}{\sqrt{Var(X_1(b_n))}} \Rightarrow N(0,\, 1)
\end{eqnarray*}
for $0<\alpha <1$. That is,  (\ref{han15}) is true.

Similarly, by (\ref{han4}) and (\ref{han6}) we have
  \begin{eqnarray*}
    \alpha- \hat{\alpha}'&=&\frac{1}{\ln b_n}\ln \Big(1+(1+o(1))\frac{2(\alpha -\hat{\alpha}')}{\hat{\alpha}'(2-\alpha)}+(1+o(1))\frac{\alpha(2-\hat{\alpha}')}{\hat{\alpha}'(2-\alpha)}\frac{(\hat{\nu}^2_n-\nu^2_n)}{\nu^2_n}\Big)\\
          &=&(1+o(1))\frac{1}{\ln b_n}\frac{(\hat{\nu}^2_n-\nu^2_n)}{\nu^2_n}
\end{eqnarray*}
for $1<\alpha <2$ and large $n$. Thus, from (\ref{han14}) of the theorem 2 it follows that
\begin{eqnarray*}
   \textbf{P}\Big( |\hat{\alpha}'-\alpha| \geq \frac{2}{n^{\beta}\ln b_n\sqrt{\ln b_n}}\Big)= \textbf{P}\Big(|\frac{\hat{\nu}^2_n -\nu^2_n}{\nu^2_n}|\geq \frac{2(1+o(1))}{n^{\beta}\sqrt{\ln b_n}}\Big)\leq \frac{2}{n^2}
 \end{eqnarray*}
for $1<\alpha <2$. By $(\alpha- \hat{\alpha}')\nu^2_n \ln b_n =(1+o(1))(\hat{\nu}^2_n-\nu^2_n)$ and (\ref{han14}), we can similarly obtain that
\begin{eqnarray*}
      \frac{\sqrt{\alpha(4-\alpha)}}{2-\alpha}\frac{\sqrt{n}(\hat{\alpha}'-\alpha)\ln b_n}{b_n^{\alpha/2}} (1+o(1)) &=&  \frac{\sqrt{n} (\alpha- \hat{\alpha}')\nu^2_n \ln b_n} {\sqrt{Var(X^2_1(b_n))}}\\
      &=&\frac{\sqrt{n}(1+o(1))(\hat{\nu}^2_n-\nu^2_n)}{\sqrt{Var(X^2_1(b_n))}} \Rightarrow N(0,\, 1)
\end{eqnarray*}
for $1<\alpha <2$. This proves (\ref{han16}).

Let $\alpha =1$ and $b_n \leq n^{1-2\beta}/\ln n$. Note that $Var(X_1(b_n)=b_n-1-(\ln b_n)^2$.  By (\ref{han7}) and the Bernstein inequality, we can get that
\begin{eqnarray*}
     \textbf{P}\Big( |\hat{\alpha}-1| \geq \frac{2\sqrt{2}}{n^{\beta}\ln b_n}\Big)&=&  \textbf{P}\Big( |n(\hat{\mu}_n- \mu_n )| \geq 2\sqrt{2}n\Big)\\
      &\leq &2\exp\{-\frac{4n}{Var(X_1(b_n))+2b_n\sqrt{2}/3}\} \leq \frac{2}{n^2}
\end{eqnarray*}
for large $n$. It follows from the Lyapunov central limit theorem that
\begin{eqnarray*}
\frac{\sqrt{n}(\hat{\alpha}-1)\ln b_n}{\sqrt{b_n}} =(1+o(1))\frac{\sqrt{n}(\hat{\mu}_n- \ln b_n )}{\sqrt{Var(X_1(b_n))}}\Rightarrow N(0,\, 1).
\end{eqnarray*}
This proves (\ref{han17}).

Let $\alpha =2$ and $b^2_n \leq n^{1-\beta}/\ln n$. By  (\ref{han8}), the Bernstein inequality and the Lyapunov central limit theorem, we can similarly prove (\ref{han18}) since $ Var(X^2_1(b_n))=(1+o(1))b_n^2$ for large $n$. It completes the proof.

\end{document}